\documentclass[twocolumn]{autart}    

\usepackage{graphicx}          

\usepackage[latin9]{inputenc}
\usepackage{amsfonts}
\usepackage{amsmath}
\usepackage{amssymb}

\begin{document}

\begin{frontmatter}

\title{Flat Sampled Data Systems, an Algorithmic Approach}

\author[Linz]{Kurt Schlacher}\ead{kurt.schlacher@jku.at},    

\address[Linz]{Institute of Automatic Control and Control Systems Technology,
               Johannes Kepler University, Linz, Austria
}  


\begin{abstract}                          
Flatness of sampled data systems can be characterized by a simple property.
They must admit the transformation to special representations, which are
the series or partial series connection of a Brunovsky normal form and a
complement. It is shown that this property follows from an integrability
condition, which must be met by the flat parametrization. The series
connections admit a simplification by reduction of the original problem
to a simpler one, which allows us to develop two algorithms, where
the first version delivers the flat outputs directly, but one has
to solve linear PDEs or nonlinear ODEs. But one gets the flat outputs
directly. The second version overcomes this problem, but one can only
test the existence of flat outputs. 
\end{abstract}

\end{frontmatter}

\section{\label{sec:Int}Introduction}
Flatness for lumped parameter systems has been introduced about 25
years ago, see e.g. \cite{FLMR} and the citations therein. It became
very popular in the control community and is an indispensable tool
today. In \cite{R1} and \cite{S1} the authors present necessary
conditions for flatness of lumped parameter time continuous systems.
Necessary and sufficient conditions are given in Levine\cite{L1},
but these conditions are not always straightforward to apply.
A constructive approach for Pfaffian systems can be found in \cite{SlSo:1},
the counterpart for explicit systems based on vector fields is shown
in \cite{SlSo-2}.

Flatness of sampled data systems can be defined analogously to the
continuous-time case, but one has the alternatives of forward shifts
or backward shifts to replace the time derivatives. The forward shift
the commonly accepted choice. The problem of input to state linearization
by static feedback has been tackeld first, see e.g. \cite{J}, \cite{G1},
\cite{AKM}, but also \cite{Ul} with all citations therein.
Extensions to exogenous presented in \cite{AM}, a more algorithmic
approach can be found in \cite{KASS}. 

This contribution uses the concept of manifolds, bundles, like tangent
or cotangent bundle, distributions, etc. Therefore, we recall
the corresponding notation and summarize some facts in Section \ref{sec:pr},
see e.g. \cite{GMS} for further details.
In Section \ref{sec:uf} we discuss two representations of sampled
data system, which are crucial for the property of flatness. One result
of this contribution is a simple Theorem, presented in Section \ref{sec:fs},
which connects integrability conditions in the space of the flat outputs
and their shifts with the representations of Section \ref{sec:uf}.
In Section \ref{sec:tst} necessary and sufficient conditions are
presented, which allow to test whether a system can be transformed
to the representation of Section \ref{sec:uf}. The series connections
admit a simplification by reduction of the original problem to a simple
one. In Section \ref{sec:alg} two algorithms are developed, where
the first one allows the determination of the flat outputs, if they
exist. The disadvantage is, that one has to solve linear PDEs or nonlinear
ODEs. Therefore, a second version is presented, where at least a test
for flatness is possible without the disadvantage of the first one.
A spin off is the fact, that flat outputs are functions of the states.
The only exception is, if the original system has redundant input.
It is worth mentioning that the first algorithm is the time discrete
counterpart to the algorithm in \cite{SlSo-2}.

\section{\label{sec:pr}Preliminaries and some technical Remarks}

In this contribution we use the geometric language of manifolds and
bundles. Let $\mathcal{M}$ denote an $m$-dimensional manifold with
local coordinates $\left(z^{1},\ldots,z^{m}\right)$, then $\mathcal{T}\left(\mathcal{M}\right)$,
$\mathcal{T}^{\ast}\left(\mathcal{M}\right)$ denote its tangent,
cotangent bundle with coordinates $\left(z^{1},\ldots,z^{m},\dot{z}^{1},\ldots,\dot{z}^{m}\right)$,
$\left(z^{1},\ldots,z^{m},\right.$ $\left. \dot{z}_{1},\ldots,\dot{z}_{m}\right)$
and canonical\footnote{Holonomic bases are constructed by the choice of
  $m$ functionally independent functions $\left(g^{1}\left(z\right),\ldots,g^{m}\left(z\right)\right)$.
  The canonical (holonomic) bases are constructed by the special choice$\left(z^{1},\ldots,z^{m}\right)$.}
bases $\left\{ \partial_{z^{1}},\ldots,\partial_{z^{m}}\right\} $,
$\left\{ \mathrm{d}z^{1},\ldots,z^{m}\right\} $. Let $C^{\infty}\left(\mathcal{M}\right)$
be the set of smooth function $f:\mathcal{M}\rightarrow\mathbb{R}$.
A smooth vector (covector) field is a smooth map $v:\mathcal{M}\rightarrow\mathcal{T}\left(\mathcal{M}\right)$,
($\omega:\mathcal{M}\rightarrow\mathcal{T}^{\ast}\left(\mathcal{M}\right)$)
or in coordinates $\dot{z}^{i}=v^{i}\left(z\right)\in C^{\infty}\left(\mathcal{M}\right)$
($\dot{z}_{i}=\omega_{i}\left(z\right)\in C^{\infty}\left(\mathcal{M}\right)$).
The canonical product $\mathcal{T}\left(\mathcal{Z}\right)\times\mathcal{T}^{\ast}\left(\mathcal{M}\right)\rightarrow C^{\infty}\left(\mathcal{M}\right)$
is denoted by\footnote{Here, we use the Einstein notation.} $v\rfloor\omega=v^{i}\omega_{i}$.
The set of all smooth vector (covector) fields is denoted by $\Gamma\left(\mathcal{T}\left(\mathcal{Z}\right)\right)$
($\Gamma\left(\mathcal{T}^{\ast}\left(\mathcal{Z}\right)\right)$).
Different coordinate systems like $z$, $\tilde{z}$ will be used
for $\mathcal{M}$, where the change is described by a diffeomorphism
$\varphi:\mathcal{M}\rightarrow\mathcal{M}$ with $\tilde{z}=\varphi\left(z\right)$. 

A distribution is a subspace $D$ of the tangent bundle $\mathcal{T}(\mathcal{M})$.
We assume, $D$ has constant rank in the neighborhood of points $z$
under consideration. A set of vector fields $B=\left\{ v_{1},\ldots,v_{n}\right\} $
is called a generator of $D$, iff $D=\mathrm{span}\left(B\right)$
is met. If the number $n$ is minimal, then $B=B_{D}$ is a basis.
A distribution $D$ is called involutive, iff it meets $\left[D,D\right]\subset D$,
where $\left[D,D\right]$ denotes the (also set valued) Lie bracket.
An involutive distribution admits a basis formed by unit vectors in
certain coordinates. Unfortunately, one has to solve linear PDEs or
nonlinear ODEs to find such a basis. But a basis $B_{D}$, which meets
$\left[v_{i},v_{j}\right]=0$, $v_{i},v_{j}\in B_{D}$ can be constructed
in a straightforward manner. Such bases, we call them adjusted, are
preferable for many calculations. 

The maximal set of symmetries $v_{i}\in\Gamma\left(\mathcal{T}\left(\mathcal{Z}\right)\right)$
of $D$ with $\left[v_{i},D\right]\subset D$ is denoted by $S\left(D\right)$.
The (Cauchy) characteristic distribution $C\left(D\right)$ is given
by $C\left(D\right)=\mathrm{span}\left(S\left(D\right)\right)\cap D$,
where $C\left(D\right)$ is involutive. A useful property is, let
$D$, $E\subset D$ be involutive distributions, then there exists
an involutive distribution $E_{c}$ such that $D=E\oplus E_{c}$ is
met.

The annihilator $D^{\perp}$ is the subset of all $\omega\in\mathcal{T^{\ast}}\left(\mathcal{Z}\right)$,
which meet $D\rfloor\omega=0$. Because of the constant rank assumption,
also $D^{\perp}$ has constant rank with $\dim\left(D^{\perp}\right)=m-n$,
$n=\dim\left(D\right)$. Iff $D$ is involutive, then $D^{\perp}$
admits a special basis $\left\{ \mathrm{\mathrm{d}}f^{m-n+1},\ldots,\mathrm{\mathrm{d}}f^{n}\right\} $
with some functions $f^{i}$, where $\mathrm{d}$ denotes the exterior
derivative on $\bigwedge\left(\mathcal{T}^{\ast}\left(\mathcal{M}\right)\right)$.
In this case $D^{\perp}$ is said to be integrable and meets $\mathrm{d}D^{\perp}\subset\mathcal{T^{\ast}}\left(\mathcal{Z}\right)\wedge D^{\perp}$.
Given two manifolds $\mathcal{M}$, $\mathcal{N}$ we consider a smooth
map $f:\mathcal{N}\rightarrow\mathcal{M}$. The pull pack $f^{\ast}\left(\omega\right)$
with $f^{\ast}:\bigwedge\left(\mathcal{T}^{\ast}\left(\mathcal{M}\right)\right)\rightarrow\bigwedge\left(\mathcal{T}^{\ast}\left(\mathcal{N}\right)\right)$
of objects $\Omega$ of the exterior algebra on $\mathcal{M}$ is
always well defined.

Let $x\in\mathbb{R}$ be a real variable, to which we assign a sequence
of values $x\left(i\right)\in\mathbb{R},$ $i=0,1,\ldots$ .The $k$-times
shift of a real variable $x$ is denoted by $x_{k}$, the assignment
$x=x\left(i\right)\in\mathbb{R}$ implies $x_{k}=x\left(i+k\right)=\sigma^{k}\left(x\left(i\right)\right)$,
$i=0,1,\ldots$ with the shift operator $\sigma$. To model a time
invariant sampled data system
\begin{eqnarray}
x_{1} & = & f\left(x,u\right)\label{eq:pr-01}
\end{eqnarray}
we introduce the bundle $\mathcal{E}\overset{\pi}{\rightarrow}\mathcal{X}$
with the $n$-dimensional base space or state manifold $\mathcal{X}$ with coordinates
$\left(x^{1},\ldots x^{n}\right)$, and the $n+m$-dimensional total
space with fiber coordinates $\left(x^{1},\ldots x^{n},u^{1},\ldots,u^{m}\right)$
and a surjective submersion $\pi:\mathcal{E}\rightarrow\mathcal{X}$.
We assume that coordinate changes for bundles respect the bundle structure
or $\tilde{x}=\varphi_{x}\left(x\right)$, $\tilde{u}=\varphi_{u}\left(x,u\right)$
with a diffeomorhism $\varphi$ is met. In this geometric picture
$f$ is a map of the type $f:\mathcal{E}\rightarrow\mathcal{X}_{1}$
with the isomorphic bundle $\mathcal{E}_{1}\overset{\pi_{1}}{\rightarrow}\mathcal{X}_{1}$.
Obviously, the shift operator $\sigma$ can be extended to geometric
objects $o$ in a straightforward manner by $\sigma\left(o\left(x,u\right)\right)=o\left(x_{1},u_{1}\right)$
for all objects defined on $\mathcal{E}$. The input distribution
$U=\mathrm{span}\left(B_{U}\right)\subset\mathcal{T}\left(\mathcal{E}\right)$,
$B_{U}=\left\{ \partial_{u^{1}},\ldots,\partial_{u^{m}}\right\} $
of the system (\ref{eq:pr-01}) meets $\pi_{\ast}$$\left(U\right)=0$
and is involutive. To avoid mathematical subtleties, we assume that
all distributions related to (\ref{eq:pr-01}) have constant rank
in the open neighborhood $\mathcal{N}$ of any point $\left(x,u\right)$,
where we develop our contribution. 

\section{\label{sec:uf}Some Useful Forms}

We assume that $f$ of (\ref{eq:pr-01}) is a surjective submersion
on $\mathcal{N}$. This is no restriction at all, otherwise the system
would not be locally reachable. Some considerations simplify, if the
system (\ref{eq:pr-01}) is transformed to the simpler form 
\begin{equation}
\begin{array}{rcl}
x_{a,1} & \;=\; & \tilde{f_{a}}\left(x,\tilde{u}\right)\\
x_{b,1} & \;=\; & \tilde{u}\;.
\end{array}\label{eq:uf-01}
\end{equation}
Possible after renumbering the equations of (\ref{eq:pr-01}) we rewrite
them as
\begin{equation}
\begin{array}{rcl}
x_{a,1} & \;=\; & f_{a}\left(x,u\right)\\
x_{b,1} & \;=\; & f_{b}\left(x,u\right)\;,
\end{array}\label{eq:uf-02}
\end{equation}
where $\mathrm{rank}\left(\partial_{u}f\right)=\mathrm{rank}\left(\partial_{u}f_{b}\right)=\tilde{m}\leq m$
is met. In the case $\tilde{m}<m$ we have $m-\tilde{m}$ redundant
inputs. The input transformation 
\begin{eqnarray*}
\tilde{u} & = & f_{b}\left(x,u\right)
\end{eqnarray*}
allows us to eliminate $u$ from $f_{a}$ of (\ref{eq:uf-02}) to
derive $\tilde{f_{a}}$ of (\ref{eq:uf-01}). If we split $u$ into
$u_{b}$, $u_{c}$ such that $\partial_{u_{b}}f_{b}\left(x,u_{b},u_{c}\right)$
is invertible to $u_{b}$, then we can assign any value to $u_{c}$. From now
on we assume that redundant inputs are eliminated, if not otherwise
mentioned.

Sometimes it is beneficial, to write $f$ of (\ref{eq:pr-01}) as
a composition of an invertible map $h$ and another submersion $g$.
From
\begin{eqnarray*}
x_{1} & = & f\left(x,u\right)\\
x_{1} & = & h\circ g\left(x,u\right)\\
h^{-1}\left(x_{1}\right) & = & g\left(x,u\right)
\end{eqnarray*}
one derives the system
\begin{eqnarray}
\tilde{x}_{1} & = & g\left(h\left(\tilde{x}\right),u\right)\label{eq:uf-03}
\end{eqnarray}
by help of the coordinate transform 
\begin{eqnarray}
x & = & h\left(\tilde{x}\right)\;.\label{eq:uf-04}
\end{eqnarray}
It is worth mentioning, that the determination of $h$ requires simple
elimination, only. Using the relations $z=g\left(x,u\right)$ to eliminate
$\left(x,u\right)$ from $f$, we obtain $h\left(z\right)$. 

An appealing form is given by
\begin{equation}
\begin{array}{rclcc}
\tilde{x}_{a,1} & = & \tilde{f_{a}}\left(\tilde{x}_{a},\tilde{x}_{b}\right)\\
\tilde{x}_{b,1} & = & \hat{f}_{b}\left(\tilde{x},u\right) & = & \tilde{u}\;,
\end{array}\label{eq:uf-05}
\end{equation}
which is the series connection of a Brunovsky normal form and a complement.
Summarizing these observations we get the following Lemma. 
\begin{lem}
\label{lem:uf-01}The system (\ref{eq:uf-02}) is transformable to
(\ref{eq:uf-05}), iff there exists a map\footnote{$\mathrm{id}_{b}$ denotes the map $\mathrm{id}_{b}\left(x_{a},x_{b}\right)=x_{b}$.}
$h=\left(h_{a},\mathrm{id}_{b}\right)$ and $\left(n-m\right)$ functions
$g$, such that 
\begin{eqnarray}
f_{a}\left(x,u\right) & = & h_{a}\left(g\left(x\right),f_{b}\left(x,u\right)\right)\label{eq:uf-06}
\end{eqnarray}
is met. In addition, the relations
\[
\begin{array}{rclcl}
\tilde{x}_{a,1} & = & g\left(h\left(\tilde{x}\right)\right) & = & \tilde{g}\left(\tilde{x}\right)\\
\tilde{x}_{b,1} & = & \hat{f}_{b}\left(h\left(\tilde{x}\right),u\right) & = & \tilde{u}
\end{array}
\]
are fulfilled.
\end{lem}
The model (\ref{eq:uf-05}) is very restrictive. Let us assume that
the system (\ref{eq:uf-02}) can be rewritten as 
\begin{equation}
\begin{array}{rclcl}
\tilde{x}_{a,1} & = & \hat{f}_{a}\left(\tilde{x},u\right) & = & \tilde{f}_{a}\left(\tilde{x},\tilde{v}\right)\\
\tilde{x}_{b_{v},1} & = & \hat{f}_{b_{v}}\left(\tilde{x},u\right) & = & \tilde{v}\\
x_{b_{u},1} & = & \hat{f}_{b_{u}}\left(\tilde{x},u\right) & = & \tilde{u}\;,
\end{array}\label{eq:uf-07}
\end{equation}
where $\mathrm{rank}\left(\partial_{u}\hat{f}_{b_{u}}\right)=m_{u}=\dim\left(\tilde{u}\right)$,
$\mathrm{rank}\left(\partial_{u}\hat{f}_{b_{v}}\right)=m_{v}=\dim\left(\tilde{v}\right)$,
$m_{v}+m_{u}=m$ and $\mathrm{rank}\left(\partial_{g}\left(\hat{f}_{a},\hat{f}_{b_{v}}\right)\right)=m_{v}$
are met. We derive a generalization of (\ref{eq:uf-05}), since by
help of the dynamic extension
\begin{eqnarray}
\tilde{v}_{1} & = & \tilde{u}_{c}\;,\label{eq:uf-08}
\end{eqnarray}
the models (\ref{eq:uf-07},\ref{eq:uf-08}) are the series connection
of a Brunovski normal form and a complement. To avoid subtleties,
we assume that $\dim\left(\tilde{v}\right)=m_{v}$ is minimal or $\dim\left(\tilde{u}\right)=m_{u}$
is maximal. Analogously to Lemma \ref{lem:uf-01} we get the following
result.
\begin{lem}
\label{lem:uf-02}The system (\ref{eq:uf-02}) is transformable to
(\ref{eq:uf-07}), iff there exists a map $h=\left(h_{a},\mathrm{id}_{b_{v}},\mathrm{id}_{b_{u}}\right)$
and $\left(n-m_{u}\right)$ functions $g=\left(g_{a},g_{v}\right)$,
$\mathrm{rank}\left(\partial_{u}g\right)=\mathrm{rank}\left(\partial_{u}g_{v}\right)=m_{v}$,
such that 
\begin{eqnarray}
f_{a}\left(x,u\right) & = & h_{a}\left(g\left(x,u\right),f_{b_{u}}\left(x,u\right)\right)\label{eq:uf-06-1}
\end{eqnarray}
is met. In addition, the relations
\[
\begin{array}{rclcl}
\tilde{x}_{a,1} & = & g_{a}\left(h\left(\tilde{x}\right),u\right) & = & \tilde{g}_{a}\left(\tilde{x},\tilde{v}\right)\\
\tilde{x}_{b_{v},1} & = & g_{b_{v}}\left(h\left(\tilde{x}\right),u\right) & = & \tilde{v}\\
x_{b_{u},1} & = & \hat{f}_{b_{u}}\left(h\left(\tilde{x}\right),u\right) & = & \tilde{u}
\end{array}
\]
are fulfilled.
\end{lem}

\section{\label{sec:fs}Properties of Flat Systems}

Let $y^{i}$ be a variable, then $y_{\left[r_{i}\right]}^{i}$ denotes
the sequence of variables
\[
y_{\left[r_{i}\right]}^{i}=\left(y_{0}^{i},y_{1}^{i},\ldots,y_{r_{i}}^{i}\right)\;,\quad y_{0}=y\;.
\]
With the sequence $\left[I\right]=\left[r_{1},\ldots,r_{m}\right]$
of non negative integers, we construct the sequences $y_{\left[I\right]}$,
$y_{\left[I_{1}\right]}$, $y_{\left[I_{-1}\right]}$, 
\begin{eqnarray*}
y_{\left[I\right]} & = & \left(y_{0}^{1},\ldots,y_{r_{i}}^{1},\ldots,y_{0}^{m},\ldots,y_{r_{m}}^{m}\right)\\
y_{\left[I_{-1}\right]} & = & \left(y_{0}^{1},\ldots,y_{r_{i}-1}^{1},\ldots,y_{0}^{m},\ldots,y_{r_{m}-1}^{m}\right)\\
y_{\left[I_{1}\right]} & = & \left(y_{1}^{1},\ldots,y_{r_{i}}^{1},\ldots,y_{1}^{m},\ldots,y_{r_{m}}^{m}\right)
\end{eqnarray*}
with $y_{\left[I\right]}\in\mathbb{R}^{R}=\mathcal{Y}_{\left[I\right]}$,
$R=\sum_{i=1}^{m}r_{i}$. The head $y_{h}$ and tail $y_{t}$ of $y_{\left[I\right]}$
are defined by
\begin{eqnarray*}
y_{h} & = & \left(y^{1},\ldots y^{m}\right)\\
y_{t} & = & \left(y_{r_{1}}^{1},\ldots y_{r_{m}}^{m}\right)\:.
\end{eqnarray*}

\begin{defn}
\label{def:fl-01}The system (\ref{eq:pr-01}) is said to be flat
(with respect to forward shifts), iff there exists a surjective submersion
$H:\mathcal{Y}_{\left[I\right]}\rightarrow\mathcal{E}$, 
\begin{equation}
\left(x,u\right)=\left(H_{x}\left(y_{\left[I\right]}\right),H_{u}\left(y_{\left[I\right]}\right)\right)=H\left(y_{\left[I\right]}\right)\,,\label{eq:fs-01}
\end{equation}
such that
\[
\sigma\left(H^{\ast}\left(x\right)\right)=H^{*}\left(f\left(x,u\right)\right)
\]
is met\footnote{The pull back $H^{\ast}\left(Z\right)$ of an indexed quantity $Z$
is a shortcut for the pullback of its elements}. The coordinates $y$ are called the coordinates of the flat outputs.
\end{defn}
\begin{rem}
\label{rem:fl-01}The relation (\ref{eq:fs-01}) is too general. We
call the map $H$ non redundant, iff $H$ is minimal with respect
to the number of coordinates and different sequences $y_{i}$ assigned
to the flat outputs generate always different sequences $\left(x_{i},u_{i}\right)$
of state $x$ and input $u$. We limit our considerations to this
type of map.
\end{rem}
Redundant inputs like $u_{c}$ of (\ref{eq:pr-01}, \ref{eq:uf-02})
are possible candidates for flat outputs. A trivial result is given
in the following remark.
\begin{rem}
\label{rem:fl-02}If the system (\ref{eq:pr-01}) is flat, then the
redundant inputs are flat outputs.
\end{rem}
By help of (\ref{eq:uf-01}) we get the relations
\begin{eqnarray*}
\sigma\left(H_{x,a}^{\ast}\left(x_{a}\right)\right) & = & f_{a}\left(H_{x,a}^{\ast}\left(x_{a}\right),H_{x,b}^{\ast}\left(x_{b}\right),\sigma\left(H_{x,a}^{\ast}\left(x_{b,1}\right)\right)\right)\\
\sigma\left(H_{x,b}^{\ast}\left(x_{a}\right)\right) & = & H_{u}^{\ast}\left(u\right)\;.
\end{eqnarray*}
It it straightforward to derive three facts:
\begin{enumerate}
\item The map $H_{x}$ meets $H_{x}=H_{x}\left(y_{\left[I_{-1}\right]}\right)$.
Since $x_{a}=H_{x,a}\left(y_{\left[I_{-1}\right]}\right)$ is the
only relation with $y=y_{0}$, $\mathrm{rank}\left(\partial_{y_{h}}H_{x}\right)=m$
must be met. Otherwise, $H$ is redundant\footnote{No shifts of any parts of $H_{x}$ can increase the $\mathrm{rank}$.
Therefore, this condition is necessary that $y$ can be expressed
as a function of $x$, $u$ and their shifts.}.
\item The map $H_{u}$ meets $H_{u}=H_{u}\left(y_{\left[I_{1}\right]}\right)$.
$H_{u}$ must depend on $y_{t}$, or $\mathrm{rank}\left(\partial_{y_{t}}H_{u}\right)\ge1$
must be met\footnote{Shifts of parts of $H_{u}$ may increase the $\mathrm{rank}$.}.
Otherwise, $H$ is redundant. 
\item The functions $H^{\ast}\left(f\left(x,u\right)\right)$ are independent
of $y_{h}$, or
\begin{eqnarray}
\partial_{y_{h}}H^{\ast}\left(f\left(x,u\right)\right) & = & 0\label{eq:fs-02}
\end{eqnarray}
is met. 
\end{enumerate}
The crucial point is Fact 3, where we have to find conditions for
the system (\ref{eq:pr-01}) such the relation (\ref{eq:fs-02}) can be met.
According to Fact 1 we introduce the spaces
\begin{eqnarray}
\mathrm{span}\left(\left\{ \mathrm{d}H_{x}^{1},\ldots,\mathrm{d}H_{x}^{n}\right\} \right) & = & X\nonumber \\
X & = & X_{y_{h}}\oplus X_{c}\;,\label{eq:fs-03}
\end{eqnarray}
where $\dim\left(X_{y_{h}}\right)=m$, $\dim\left(X_{c}\right)=n-m$
and $\partial_{y^{i}}\rfloor X_{c}=\mathrm{\mathrm{span}\left(\left\{ 0\right\} \right)}$,
$i=1,\ldots,m$ and are met. Since $X_{c}$ is maximal with respect
to the dimension, we derive the following Lemma.
\begin{lem}
\label{lem:fs-01}The system (\ref{eq:pr-01}) is (locally) transformable
to (\ref{eq:uf-05}), iff the space $X_{c}$ of (\ref{eq:fs-03})
is integrable.
\end{lem}
Integrability of $X_{c}$ implies there exists a basis of the form
$\left\{ \mathrm{d}G^{1}\left(H_{x}\right),\ldots,\mathrm{d}G^{n-m}\left(H_{x}\right)\right\} $,
where $G^{i}$ is independent of $y_{h}$. Because of Definition \ref{def:fl-01},
there exists $n-m$ function $g^{i}\left(x\right)$, such that $H^{\ast}\left(g^{i}\right)=G^{i}$
is meet. By help of Lemma \ref{lem:uf-01} we are done. Obviously,
the relation (\ref{eq:fs-02}) is met.

If $X_{c}$ is not integrable, one can try to augment $X_{c}$ by
adding a subspace of
\begin{eqnarray*}
U & = & \mathrm{span}\left(\left\{ \mathrm{d}H_{u}^{1},\ldots,\mathrm{d}H_{u}^{m}\right\} \right)\;.
\end{eqnarray*}
Let us consider the case $\dim\left(U\right)=1$ first. Since $\partial_{y_{t}^{1}}\rfloor\mathrm{d}H_{u}^{1}\neq0$
is met because of Fact 2, the space $X_{c}\oplus U$ is integrable,
iff $X_{c}$ is integrable. We gain nothing and cannot
meet relation (\ref{eq:fs-02}). This is stated by the following
Lemma.
\begin{lem}
\label{lem:fs-02}A necessary condition for the system (\ref{eq:uf-01})
with $\dim\left(u\right)=1$ to be flat is, that it is (locally) transformable
to the representation (\ref{eq:uf-05}) .
\end{lem}
The case $\dim\left(U\right)=m>1$ is similar to the above one, iff
$m=m_{u}$ with $\mathrm{rank}\left(\partial_{y_{t}}H_{u}\right)=m_{u}$
is met. If $m_{u}<m$ is met, we can split $U$, possible after
renumbering the functions $H_{u}$ in the following manner
\begin{eqnarray*}
U & = & U_{y_{t}}\oplus U_{c}\\
U_{y_{t}} & = & \mathrm{span}\left(\left\{ \mathrm{d}H_{u}^{1},\ldots,\mathrm{d}H_{u}^{m_{u}}\right\} \right)\;,
\end{eqnarray*}
such that $\mathrm{rank}\left(\partial_{y_{t}}\left(H_{u}^{1},\ldots,H_{u}^{m_{u}}\right)\right)=m_{u}$ and  $\partial_{y_{t}^{i}}\rfloor U_{c}=\mathrm{span}\left(\left\{ 0\right\} \right)$
are fulfilled. Now $X_{c}\oplus U_{c}$ must be integrable, or there exists
a basis of the form $\left\{ \mathrm{d}G^{1}\left(H_{x},H_{u}\right),\ldots,\right. $
$\left.\mathrm{d}G^{n-m_{u}}\left(H_{x},H_{u}\right)\right\} $,
where $G^{i}$ is independent of $y_{h}$ . Because of Definition
\ref{def:fl-01}, there exists $n-m_{u}$ function $g^{i}\left(x,u\right)$,
such that $H^{\ast}\left(g^{i}\right)=G^{i}$ is met. By construction
the functions\textbf{ $G^{i}$} are independent of $y_{t}$, too,
or $\partial_{H_{u}}G\partial_{y_{t}}H_{u}=0$ with $\mathrm{rank}\left(\partial_{y_{t}}H_{u}\right)=m_{u}$
is met. Since $H_{u}$ is a submersion the relation $\partial_{H_{u}}G=\partial_{u}g\left(x,u\right)_{u=H_{u}}$,
implies $\mathrm{rank}\left(\partial_{H_{u}}G\right)=\mathrm{rank}\left(\partial_{u}g\right)=m-m_{u}$.
Finally, by help of Lemma \ref{lem:uf-02}, we derive the following
theorem.
\begin{thm}
\label{thm:fs-01}A necessary condition for the system (\ref{eq:pr-01})
to be flat is, that it is (locally) transformable to the representation
(\ref{eq:uf-07}) .
\end{thm}

\section{\label{sec:tst}Differential Geometry Approaches }

In Section \ref{sec:uf} two system representations have been discussed,
which are a series or partial series connection of a Brunovski normal
form and a complement. Now we add the missing tests to check whether
a transformation to these forms is possible. Let us consider the involutive
distribution $K$, 
\begin{equation}
K=\left\{ v\in\Gamma\left(\mathcal{T}\left(\mathcal{E}\right)\right)|v\left(f^{i}\right)=0,\:i=1,\ldots,n\right\} \label{eq:tst-01}
\end{equation}
for (\ref{eq:pr-01}). 

Let $\tilde{K}$ denote this distribution in the coordinates of (\ref{eq:uf-05}).
Obviously, $\tilde{K}$ has a basis, which is independent of $\tilde{u}$.
This fact implies $\left[\tilde{U},\tilde{K}\right]\subseteq\tilde{U}\oplus\tilde{K}$
or 
\begin{eqnarray}
\left[U,K\right] & \subseteq & U\oplus K\label{eq:tst-02}
\end{eqnarray}
in the coordinates of (\ref{eq:pr-01}). The relation (\ref{eq:tst-02})
implies that the annihilator $\left(U\oplus K\right)^{\perp}$ has
a basis of differentials of $\left(n-m\right)$ functions $g\left(x\right)$,
see Lemma \ref{lem:uf-01}, and we are done. Summarizing we get the
following Lemma.
\begin{lem}
\label{lem:tst-01}The relation (\ref{eq:tst-02}) is necessary and
locally sufficient that the system (\ref{eq:pr-01}) is transformable
to (\ref{eq:uf-05}). 
\end{lem}
To derive a test for a transformation to (\ref{eq:uf-07}), we simple
state that $\tilde{K}$ has a basis, which is independent of $\tilde{u}$.
This can be expressed as $\left[\tilde{U}_{b_{u}},\tilde{K}\right]\subseteq\tilde{U}_{b_{u}}\oplus\tilde{K}$
for the involutive subdistribution $\tilde{U}_{b_{u}}\subset\tilde{U}$.
In the coordinates of (\ref{eq:pr-01}) we rewrite this condition
as
\begin{eqnarray}
\left[U_{b_{u}},U_{b_{u}}\oplus K\right] & \subseteq & U_{b_{u}}\oplus K\label{eq:tst-03}
\end{eqnarray}
for $U_{b_{u}}\subset U$, $\dim\left(U_{b_{u}}\right)=m_{u}$. The
relation (\ref{eq:tst-03}) implies that the annihilator $\left(U_{b_{u}}\oplus K\right)^{\perp}$
has a basis of differentials of $\left(n-m_{u}\right)$ functions
$g=\left(g_{a}\left(x,u\right),g_{v}\left(x,u\right)\right)$ with
$\mathrm{rank}\left(\partial_{u}g\right)=\mathrm{rank}\left(\partial_{u}g_{v}\right)=m_{v}=m-m_{u}$.
By use of Lemma \ref{lem:uf-02} we get the following result.
\begin{lem}
\label{lem:tst-02}The relation (\ref{eq:tst-03}) is necessary and
locally sufficient that the system (\ref{eq:pr-01}) is transformable
to (\ref{eq:uf-07}). 
\end{lem}
The system (\ref{eq:uf-05}) is a serial connection of a Brunovski
normal form $B$ and a complement $C$. The distribution of inputs
of $B$ generated by $\tilde{B}_{B}=\left\{ \partial_{\tilde{u}^{1}},\ldots,\partial_{\tilde{u}^{m}}\right\} $
and of $C$ generated by $\tilde{B}_{\tilde{X}_{b}}=\left\{ \partial_{\tilde{x}_{b}^{1}},\ldots,\partial_{\tilde{x}_{b}^{m}}\right\} $
are connected by the push forward $\tilde{f}_{\ast}$ given by
\begin{eqnarray*}
\dot{\tilde{x}}_{1} & = & \left[\begin{array}{cc}
X & 0\\
Y & I
\end{array}\right]\left[\begin{array}{c}
\dot{\tilde{x}}\\
\dot{\tilde{u}}
\end{array}\right]\;.
\end{eqnarray*}
This implies $\text{\ensuremath{\partial_{\tilde{x}_{b,1}^{j}}=}}\tilde{f}_{\ast}\partial_{\tilde{u}^{j}}$
and $\partial_{\tilde{x}_{b}^{j}}=\sigma^{-1}\left(\partial_{\tilde{x}_{b,1}^{j}}\right)$,
where $\sigma^{-1}$ is well defined. 

Now, we construct the equivalent relations in the coordinates of (\ref{eq:uf-02})
and start with an involutive distribution on inputs generated $B_{B}=\left\{ \eta_{1},\ldots,\eta_{m}\right\} $.
We assume $\left[\eta_{i},\eta_{j}\right]=0$ and adapt the push forward
operation to $B_{B}$ such that 
\begin{eqnarray*}
\dot{x}_{1} & = & \left[\begin{array}{cc}
X & F_{a}\\
Y & F_{b}
\end{array}\right]\left[\begin{array}{c}
\dot{x}\\
\dot{\eta}
\end{array}\right]
\end{eqnarray*}
is met. Let us assume, $F_{b}$ is invertible. By help of $\dot{\tilde{\eta}}=F_{b}\dot{\eta}$,
we adapt $B_{B}$ to $\tilde{B}_{B}=\left\{ \tilde{\eta}_{1},\ldots,\tilde{\eta}_{m}\right\} $
and get
\begin{eqnarray}
\dot{x}_{1} & = & \left[\begin{array}{cc}
X & F_{a}F_{b}^{-1}\\
Y & I
\end{array}\right]\left[\begin{array}{c}
\dot{x}\\
\dot{\tilde{\eta}}
\end{array}\right]\;.\label{eq:tst-04}
\end{eqnarray}
The field $\xi_{j}=f_{\ast}\tilde{\eta}_{j}=$ splits into to two
parts $\xi_{j,a}+\xi_{j,b}$, where $\sigma^{-1}\left(\xi_{j,b}\right)$
is well defined. From above we know that the part $\xi_{j,a}$ vanishes
in certain coordinates. Therefore, $\sigma^{-1}\left(\xi_{j}\right)$
is well defined, too, but the functions of (\ref{eq:pr-01}) are needed
to perform $\sigma^{-1}$ in general. Since the case of (\ref{eq:uf-05}) is
almost identical, it is omitted.

\section{\label{sec:alg}An Algorithm}

The Theorem \ref{thm:fs-01} allows as to construct a simple algorithm
for the determination of flat outputs, if they exist. A necessary
condition for a system $S$ like (\ref{eq:pr-01}) to be flat is,
it admits a (partial) series connection of a Brunovsky normal form
$B$ and a complement $C$, where $S=\left(B,C\right)$ is the shortcut
for the series connection. The set of all inputs of $S$ is denoted
by $A_{S}$ the set of redundant inputs is denoted by $Y_{S}$, where
$A_{S}=Y_{S}\cup U_{S}$. Its non redundant inputs $U_{S}$ split
into $U_{S}=U_{B}\cup V_{S}$, where $U_{B}$ is the input of $B$
and $V_{S}$ is its complement. The input of $C$ is given by $X_{B}\cup V_{S}$,
where $X_{B}$ denotes the set of states of $B$. The empty system
has no state equations, but may have redundant variables.
The following algorithm is based on the following observation.
\begin{lem}
\label{lem:al-01}If the system $S=\left(B,C\right)$ is flat, then
$C$ is flat, too.
\end{lem}
This Lemma follows from the following facts. If $S$ is flat, then
$Y_{S}$, $U_{B}$, $V_{S}$ and the state $X_{S}$ are uniquely determined
by the flat outputs. The states $X_{S}$ and $X_{B}\cup X_{C}$ are
connected by a simple state transformation. Therefore, the system
$C$ with input $X_{B}\cup V_{S}$ and state $X_{C}$ must be flat,
too.

Lemma \ref{lem:al-01} allows us to derive the following algorithm,
which determines a set flat outputs of a system $S$, if it is flat.
\begin{enumerate}
\item[0)] Set $Y=\left\{ \right\} $.
\item[1)] Determine $Y_{S}$, $A_{S}=U_{S}\cup Y_{S}$. Set $Y=Y\cup Y_{S}$.
\item[2)] If $S=\varnothing$ stop. The system is flat with the flat outputs
$Y$. 
\item[3)] If $S=\left(B,C\right)$ exists continue. Otherwise stop, the system
is not flat.
\item[4)] Determine $V_{S}$, $U_{B}$ with $U_{S}=U_{B}\cup V_{S}$. Set $A=X_{B}\cup V_{S},$
$S=C$, $A_{S}=A$. Goto 1.
\end{enumerate}
A simple consequence of this algorithm is the following Lemma, see
also Remark \ref{rem:fl-02}.
\begin{lem}
\label{lem:al-02}If the system $S$ is flat and meets $Y_{S}=\left\{ \right\} $,
then the flat outputs are functions of the states $X_{S}$ only.
\end{lem}
Please, note that flat outputs are related to redundant inputs. If
$Y_{S}=\left\{ \right\} $ is met, then the inputs of $C$, $S=\left(B,C\right)$,
follow as $X_{B}$. Therefore, possible redundant inputs $Y_{C}$
are functions of $X_{B}$. A simple repetition of these arguments
proves Lemma \ref{lem:al-02}.

The disadvantage of the algorithm from above is that one has to solve
linear PDEs or nonlinear ODEs to determine the series connection.
This will be illustrated in the examples. Often it is enough to check,
whether a system is flat. The following algorithm is a copy of the
previous one, where sets of variables are replaced by distributions.
E.g. $A_{S}$ denotes the set of all inputs $\left\{ u_{1},\ldots.,u_{m}\right\} $
in the previous algorithm and the set $\left\{ \partial_{u^{1}},\,\ldots,\partial_{u^{m}}\right\} $
in the algorithm below. Furthermore, we use adjusted bases in the
examples to simplify the calculations.
\begin{enumerate}
\item[0)] Set $Y=R=\mathrm{span}\left(\left\{ 0\right\} \right)$, $\kappa=\left|X_{S}\right|$.
\item[1)] Determine $Y_{S}$, $A_{S}=U_{S}\oplus Y_{S}$. Set $Y=Y\oplus Y_{S}$. 
\item[2)] If $\kappa=0$ stop. The system is flat.
\item[3)] If $\dim\left(\hat{U}\right)>0$ for $\left[\hat{U}\oplus R,K\right]\subset\hat{U}\oplus R\oplus K$,
$\hat{U}\subset U_{S}$ is met continue. Otherwise stop, the system
is not flat.
\item[4)] Determine $V_{S}$, $U_{S}=\hat{U}\oplus V_{S}$. Set $A_{S}=\Pi\left(\hat{U}\right)\oplus V$,
$R=R\oplus\hat{U}$, $\kappa=\kappa-\dim\left(\hat{U}\right)$. Goto
2.
\end{enumerate}
Several facts are worth mentioning. To derive the flat outputs, one
has to integrate $R^{\perp}$, which requires to solve linear PDEs
or nonlinear ODEs. In addition, one has to construct a suitable basis
for $Y$ to derive the flat outputs. In general, this requires the
knowledge of the generating functions of $R^{\perp}$. Furthermore,
all calculations are done with all inputs of the Brunovski forms,
which are characterized by the distribution $R$ to construct a suitable
basis. 

Two examples are considered, where the first one demonstrates which
operations are required to determine a flat output or to show its
existence only. The second one is little bit more challenging. 

\subsection{Example 1}

First we consider the simple example 
\begin{equation}
\begin{array}{rcl}
x_{1}^{1} & = & \left(x^{1}+x^{2}\right)^{3}x^{2}u\\
x_{1}^{2} & = & x^{2}u
\end{array}\label{eq:ex-01}
\end{equation}
with one input $u$ and input distribution $U=\mathrm{span}\left(\left\{ \partial_{u}\right\} \right)$
to demonstrate, which operations are required to determine a flat
output or to show its existence, only. We set $Y=\left\{ \right\} $
and start the algorithm.
\begin{enumerate}
\item $Y_{S}=\left\{ \right\} $, $U_{S}=\left\{ u\right\} $.
\item $S\ne\varnothing$.
\item The conditions of Lemma \ref{lem:uf-01} are met by
\begin{eqnarray*}
g^{1}\left(x^{1},x^{2}\right) & = & x^{1}+x^{2}\\
g^{2}\left(x^{1},x^{2},u\right) & = & f_{2}=x^{2}u\:.
\end{eqnarray*}
To derive the map $h$, see (\ref{eq:uf-04}), we solve the set of
equations 
\begin{eqnarray*}
z^{1} & = & x^{1}+x^{2}\\
z^{2} & = & x^{2}u
\end{eqnarray*}
for $u$ and some $x$, here $x^{1}$ , and get 
\begin{eqnarray*}
h^{1}\left(z\right) & = & \left(z^{1}\right)^{3}z^{2}\\
h^{2}\left(z\right) & = & z^{2}\;.
\end{eqnarray*}
The transformation is given by
\begin{eqnarray*}
x^{1} & = & \left(\tilde{x}^{1}\right)^{3}\tilde{x}^{2}\\
x^{2} & = & \tilde{x}^{2}
\end{eqnarray*}
and the transformed system reads as 
\begin{eqnarray*}
\tilde{x}_{1}^{1} & = & \left(\left(\tilde{x}^{1}\right)^{3}+1\right)\tilde{x}_{2}\\
\tilde{x}_{1}^{2} & = & \tilde{x}^{2}u\;.
\end{eqnarray*}
By help of the input transformation $\tilde{u}=\tilde{x}^{2}u$, we
get the Brunovski form
\begin{eqnarray*}
\tilde{x}_{1}^{2} & = & \tilde{u}
\end{eqnarray*}
and the complement 
\begin{eqnarray*}
\tilde{x}_{1}^{1} & = & \left(\left(\tilde{x}^{1}\right)^{3}+1\right)\tilde{x}^{2}\;.
\end{eqnarray*}
\item $V_{S}=\left\{ \right\} $, $U_{B}=\left\{ \tilde{u}\right\} $. $A_{S}=\left\{ \tilde{x}^{2}\right\} $.
\end{enumerate}
We repeat the procedure:
\begin{enumerate}
\item $Y_{S}=\left\{ \right\} $, $U_{S}=\left\{ \tilde{x}_{2}\right\} $.
\item $S\ne\varnothing$.
\item By help of the input transformation $\tilde{\tilde{x}}^{2}=\left(\left(\tilde{x}^{1}\right)^{3}+1\right)\tilde{x}^{2}$,
we get the Brunovski form
\begin{eqnarray*}
\tilde{x}_{1}^{1} & = & \tilde{\tilde{x}}^{2}
\end{eqnarray*}
and its complement is the empty system.
\item $V_{S}=\left\{ \right\} $, $U_{B}=\left\{ \tilde{\tilde{x}}^{2}\right\} $.
$A_{S}=\left\{ \tilde{x}^{1}\right\} $. 
\end{enumerate}
We repeat the procedure:
\begin{enumerate}
\item $Y_{S}=\left\{ \tilde{x}^{1}\right\} $, $U_{S}=\left\{ \right\} $,
$Y=Y_{S}$
\item $S=\varnothing$. The system is flat.
\end{enumerate}
We derive the flat output in the transformed coordinates. In the original
ones we get $y=\tilde{x}^{1}=\left(x^{1}/x^{2}\right)^{1/3}$ . 

\bigskip{}

To check only, whether the system is flat, we determine $K$,
\begin{eqnarray*}
K & = & \mathrm{span}\left(\left\{ x^{2}\partial_{x^{1}}-x^{1}\partial_{x^{2}}+u\partial_{u}\right\} \right)
\end{eqnarray*}
and set $Y=R=\mathrm{span}\left(\left\{ 0\right\} \right)$, $\kappa=2$
and start the algorithm. 
\begin{enumerate}
\item $Y_{S}=\mathrm{span}\left(\left\{ 0\right\} \right)$, $U_{S}=\mathrm{span}\left(\left\{ \partial_{u}\right\} \right)$,
$Y=Y_{S}\oplus Y$.
\item $\kappa=2\ne0$.
\item $\left[U_{S},K\right]\subset U_{S}\oplus K$ is met. 
\item $V_{S}=\left\{ \mathrm{span}\left(\left\{ 0\right\} \right)\right\} $,
$U_{B}=U_{S}$. We choose $f_{b}=f^{2}$ and derive 
\begin{eqnarray*}
\xi & = & \left(x^{1}+x^{2}\right)^{3}\partial_{x_{1}^{1}}+\partial_{x_{1}^{2}}\\
 & = & \frac{x_{1}^{1}}{x_{1}^{2}}\partial_{x_{1}^{1}}+\partial_{x_{1}^{2}}\mod(\ref{eq:ex-01})\\
\sigma^{-1}\left(\xi\right) & = & \frac{x^{1}}{x^{2}}\partial_{x^{1}}+\partial_{x^{2}}
\end{eqnarray*}
according to formula (\ref{eq:tst-04}). We set 
$R=R\oplus U_{B}$, $\kappa=\kappa-1$, $A_{S}=\mathrm{span}\left(\left\{ x^{1}\partial_{x^{1}}+x^{2}\partial_{x^{2}}\right\} \right)$. 
\end{enumerate}
Now we repeat the procedure:
\begin{enumerate}
\item $Y_{S}=\mathrm{span}\left(\left\{ 0\right\} \right)$, $U_{S}=A_{S}$,
$Y=Y_{S}\oplus Y$. 
\item $\kappa=1\ne0$.
\item $\left[U_{S},K\right]\subset U_{S}\oplus K$ is met. 
\item $V_{S}=\left\{ \mathrm{span}\left(\left\{ 0\right\} \right)\right\} $,
$U_{B}=U_{S}$. According to formula (\ref{eq:tst-04}) with $f_{b}=\left(f^{1},f^{2}\right)$
and the sequence $\left(x^{1}\partial_{x^{1}}+x^{2}\partial_{x^{2}},\partial_{u}\right)$
we get
\begin{eqnarray*}
\xi_{1} & = & \partial_{x_{1}^{1}}\\
\xi_{2} & = & \partial_{x_{1}^{2}}\;.
\end{eqnarray*}
We set $A_{S}=\mathrm{span}\left(\left\{ \partial_{x^{1}}\right\} \right)$,
$R=R\oplus U_{B}$, $\kappa=\kappa-1$. 
\end{enumerate}
Now, we repeat the procedure:
\begin{enumerate}
\item $Y_{S}=\mathrm{span}\left(\left\{ \partial_{x^{1}}\right\} \right)$,
$U_{S}=\mathrm{span}\left(\left\{ 0\right\} \right)$, $Y=Y_{S}\oplus Y$. 
\item $\kappa=0$. The system is flat.
\end{enumerate}
The annihilator of $R=\mathrm{span}\left(\left\{ \partial_{u},x^{1}\partial_{x^{1}}+x^{2}\partial_{x^{2}}\right\} \right)$
follows as $\mathrm{span}\left(\left\{ \mathrm{d}p\left(x^{1}/x^{2}\right)\right\} \right)$
with an arbitrary function $p\left(\cdot\right)$. The field $\partial_{x^{1}}$
is not a symmetry of $R$, but $x^{2}\partial_{x^{1}}$ is one. By
help of the additional equation
\begin{eqnarray*}
x^{2}\partial_{x^{1}}p\left(x^{1}/x^{2}\right) & = & 1
\end{eqnarray*}
 we get the special flat output $y=x^{1}/x^{2}$.

\subsection{Example 2}

A more interesting example is given by the system
\begin{equation}
\begin{array}{rcl}
x_{1}^{1} & = & \frac{x^{2}}{x^{4}+x^{1}\left(u^{1}-u^{2}\right)}\\
x_{1}^{2} & = & x^{4}+x^{1}\left(u^{1}-u^{2}\right)\\
x_{1}^{3} & = & x^{2}u^{1}+x^{4}\\
x_{1}^{4} & = & x^{2}u^{1}+x^{3}\left(u^{2}-u^{1}\right)+x^{4}\left(u^{1}-u^{2}+1\right)
\end{array}\label{eq:ex-02}
\end{equation}
with four states and two inputs. First, we determine $K$,
\begin{eqnarray*}
K & = & \mathrm{span}\left(\left\{ k_{1},k_{2}\right\} \right)\\
  k_{1} & = & \frac{x^{2}-x^{1}}{u^{1}-u^{2}}\partial_{x^{1}}+\frac{\left(u^{2}-u^{1}\right)x^{2}-x^{3}+x^{4}}{u^{1}-u^{2}}\partial_{x^{3}}\\
& &              +x^{2}\partial_{x^{4}}+\partial_{u^{1}}\\
k_{2} & = & \frac{x^{1}}{u^{1}-u^{2}}\partial_{x^{1}}+\frac{x^{3}-x^{4}}{u^{1}-u^{2}}\partial_{x^{3}}+\partial_{u^{2}}\;.
\end{eqnarray*}
We set $Y=\left\{ \right\} $ and start the algorithm.
\begin{enumerate}
\item $Y_{S}=\left\{ \right\} $, $U_{S}=\left\{ u^{1},u^{2}\right\} $.
\item $S\ne\varnothing$.
\item Since the requirements of Lemma \ref{lem:uf-01} are not met, we try
to meet Lemma \ref{lem:uf-02}. According to (\ref{eq:tst-03}) we
derive 
\begin{eqnarray*}
U_{b_{v}} & = & \mathrm{span}\left(\left\{ \partial_{u^{1}}+\partial_{u^{2}}\right\} \right)\;.
\end{eqnarray*}
The annihilator of $U_{b_{v}}\oplus K$ is $\mathrm{span}\left(\mathrm{d}g^{1},\mathrm{d}g^{2},\mathrm{d}g^{3}\right)$,
\begin{eqnarray*}
g^{1} & = & x^{2}\\
g^{2} & = & x^{4}+x^{1}\left(u^{1}-u^{2}\right)\\
g^{3} & = & \left(u^{1}-u^{2}\right)\left(x^{3}-x\right)\;.
\end{eqnarray*}
Next, we solve the equations 
\begin{eqnarray*}
z^{1} & = & x^{2}\\
z^{2} & = & x^{4}+x^{1}\left(u^{1}-u^{2}\right)\\
z^{3} & = & \left(u^{1}-u^{2}\right)\left(x^{3}-x^{4}\right)\\
z^{4} & = & f_{4}=x^{2}u^{1}+x^{3}\left(u^{2}-u^{1}\right)+x^{4}\left(u^{1}-u^{2}+1\right)
\end{eqnarray*}
with respect to $x^{1}$, $x^{2}$, $x^{3}$, $x^{4}$ and derive
the map $h$
\begin{eqnarray*}
h^{1}\left(z\right) & = & \frac{z^{1}}{z^{2}}\\
h^{2}\left(z\right) & = & z^{2}\\
h^{3}\left(z\right) & = & z^{3}+z^{4}\\
h^{3}\left(z\right) & = & z^{4}
\end{eqnarray*}
together with the transformation 
\begin{eqnarray*}
x^{1} & = & \frac{\tilde{x}^{1}}{\tilde{x}^{2}}\\
x^{2} & = & \tilde{x}^{2}\\
x^{3} & = & \tilde{x}^{3}+\tilde{x}^{4}\\
x^{4} & = & \tilde{x}^{4}\:.
\end{eqnarray*}
The transformed system reads as
\begin{eqnarray*}
\tilde{x}_{1}^{1} & = & \tilde{x}^{2}\\
\tilde{x}_{1}^{2} & = & \tilde{x}^{4}+\frac{\tilde{x}^{1}}{\tilde{x}^{2}}\left(u^{1}-u^{2}\right)\\
\tilde{x}_{1}^{3} & = & \tilde{x}^{3}\left(u^{1}-u^{2}\right)\\
\tilde{x}_{1}^{4} & = & \tilde{x}^{4}+\left(\tilde{x}^{2}-\tilde{x}^{3}\right)u^{1}+\tilde{x}^{2}\tilde{x}^{3}\;.
\end{eqnarray*}
The input transformation 
\begin{eqnarray*}
\tilde{u}_{1} & = & \tilde{x}^{4}+\left(\tilde{x}^{2}-\tilde{x}^{3}\right)u^{1}+\tilde{x}^{2}\tilde{x}^{3}\\
\tilde{v}_{1} & = & \tilde{x}^{3}\left(u^{1}-u^{2}\right)
\end{eqnarray*}
leads to the Brunovski form
\begin{eqnarray*}
\tilde{x}_{1}^{4} & = & \tilde{u}^{1}
\end{eqnarray*}
and the complement
\begin{eqnarray*}
\tilde{x}_{1}^{1} & = & \tilde{x}^{2}\\
\tilde{x}_{1}^{2} & = & \tilde{x}^{4}+\frac{\tilde{x}^{1}}{\tilde{x}^{2}\tilde{x}^{3}}\tilde{v}^{1}\\
\tilde{x}_{1}^{3} & = & \tilde{v}^{1}\;.
\end{eqnarray*}
 
\item $V_{S}=\left\{ \tilde{v}^{1}\right\} $, $U_{B}=\left\{ \tilde{u}^{1}\right\} $.
$A_{S}=\left\{ \tilde{v}^{1},\tilde{x}^{4}\right\} $.
\end{enumerate}
The next step would be to repeat the previous procedure. By help of
the transformation
\begin{eqnarray*}
\tilde{u}^{2} & = & \tilde{x}^{4}+\frac{\tilde{x}^{1}}{\tilde{x}^{2}\tilde{x}^{3}}\tilde{v}^{1}
\end{eqnarray*}
we rewrite the complement as 
\begin{eqnarray*}
\tilde{x}_{1}^{1} & = & \tilde{x}^{2}\\
\tilde{x}_{1}^{2} & = & \tilde{u}^{2}\\
\tilde{x}_{1}^{3} & = & \tilde{v}^{1}\;.
\end{eqnarray*}
 This is already a Brunovsky normal form and we see that $y^{1}=\tilde{x}^{1}=x^{1}x^{2}$,
$y^{2}=\tilde{x}^{3}=x^{3}-x^{4}$ are flat outputs.\bigskip{}

To check only, whether the system is flat, we set $Y=R=\mathrm{span}\left(\left\{ 0\right\} \right)$,
$\kappa=4$ and start the algorithm. 
\begin{enumerate}
\item $Y_{S}=\mathrm{span}\left(\left\{ 0\right\} \right)$, $U_{S}=\mathrm{span}\left(\left\{ \partial_{u^{1}},\partial_{u^{2}}\right\} \right)$,
$Y=Y_{S}\oplus Y$. 
\item $\kappa=4\ne0$.
\item $\left[U_{b},K\right]\subseteq U_{b}\oplus K$ is met for $U_{b}=\mathrm{\mathrm{span}}\left(\left\{ \partial_{u^{1}}+\partial_{u^{1}}\right\} \right)$. 
\item $V_{S}=\left\{ \mathrm{span}\left(\left\{ \partial_{u^{1}}\right\} \right)\right\} $,
$U_{B}=U_{b}$. According to (\ref{eq:tst-04}) we get
\begin{eqnarray*}
\xi & = & \partial_{x_{1}^{3}}+\partial_{x_{1}^{4}}
\end{eqnarray*}
for the choice $f_{b}=f^{4}$ . We set $A_{S}=V_{S}\oplus\mathrm{span}\left(\left\{ \partial_{x^{3}}+\partial_{x^{4}}\right\} \right)$,
$R=R\oplus U_{B}$, $\kappa=\kappa-1$. 
\end{enumerate}
Now we repeat the procedure.
\begin{enumerate}
\item $Y_{S}=\mathrm{span}\left(\left\{ 0\right\} \right)$, $U_{S}=A_{S}$,
$Y=Y_{S}\oplus Y$. 
\item $\kappa=3\ne0$.
\item $\left[U_{S},K\right]\subset U_{S}\oplus K$ is met. 
\item $U_{S}=\left\{ \mathrm{span}\left(\left\{ 0\right\} \right)\right\} $,
$U_{B}=U_{S}$. With the choice $f_{b}=f^{2},f^{3},f^{4}$ and the
sequence $\partial_{u^{1}},\partial_{x^{3}}+\partial_{x^{4}},\partial_{u^{1}}+\partial_{u^{2}}$
we derive
\begin{eqnarray*}
\xi_{1} & = & -\frac{x_{1}^{1}}{x_{1}^{2}}\partial_{x_{1}^{1}}+\partial_{x_{1}^{2}}\\
\xi_{2} & = & \partial_{x_{1}^{3}}\;.\\
\xi_{3} & = & \partial_{x_{1}^{4}}
\end{eqnarray*}
We set $A_{S}=\mathrm{span}\left(\left\{ x^{1}\partial_{x^{1}}-x^{2}\partial_{x^{2}},\partial_{x^{3}}\right\} \right)$,
$R=R\oplus U_{B}$, $\kappa=\kappa-2$. 
\end{enumerate}
Now we repeat the procedure.
\begin{enumerate}
\item $Y_{S}=\mathrm{span}\left(\left\{ \partial_{x^{3}}\right\} \right)$,
$U_{S}=\mathrm{span}\left(\left\{ x^{1}\partial_{x^{1}}-x^{2}\partial_{x^{2}}\right\} \right)$,
$Y=Y_{S}\oplus Y$. 
\item $\kappa=1\ne0$.
\item $\left[U_{S},K\right]\subset U_{S}\oplus K$ is met. 
\item $V_{S}=\left\{ \mathrm{span}\left(\left\{ 0\right\} \right)\right\} $,
$U_{B}=U_{S}$. With $f_{b}=f^{1},f^{2},f^{3},f^{4}$ and the sequence
$x^{1}\partial_{x^{1}}-x^{2}\partial_{x^{2}},\partial_{u^{1}},$ $\partial_{x^{3}}+\partial_{x^{4}},\partial_{u^{1}}+\partial_{u^{2}}$
we derive 
\begin{eqnarray*}
\xi & = & \partial_{x^{1}}\;.
\end{eqnarray*}
We set $A_{S}=\mathrm{span}\left(\left\{ \partial_{x^{1}}\right\} \right)$,
$R=R\oplus U_{B}$, $\kappa=\kappa-1$. 
\end{enumerate}
Now we repeat the procedure.
\begin{enumerate}
\item $Y_{S}=Y\mathrm{span}\left(\left\{ \partial_{x^{1}}\right\} \right)$,
$U_{S}=\mathrm{span}\left(\left\{ 0\right\} \right)$, $Y=Y_{S}\oplus Y$. 
\item $\kappa=0$. The system is flat.
\end{enumerate}
The annihilator of $R=\mathrm{span}\left(\left\{ x^{1}\partial_{x^{1}}-x^{2}\partial_{x^{2}},\partial_{x^{3}}+\partial_{x^{4}},\right. \right.$
    $\left. \left. \partial_{u^{1}},\partial_{u^{1}}+\partial_{u^{2}}\right\} \right)$
    follows as $\mathrm{span}\left(\left\{ \mathrm{d}p\left(x^{1}x^{2},x^{3}-x^{4}\right),\right.\right.$
        $\left.\left.\mathrm{d}p\left(x^{1}x^{2},x^{3}-x^{4}\right)\right\} \right)$
with arbitrary functions $p\left(\cdot,\cdot\right)$,  $q\left(\cdot,\cdot\right)$.
With $Y=\mathrm{span}\left(\left\{ \partial_{x^{1}},\partial_{x^{3}}\right\} \right)$
we see that $\partial_{x^{1}}$ is not a symmetry of $R$, but $1/x^{2}\partial_{x^{1}}$
is one. By help of the additional equation
\begin{eqnarray*}
\frac{1}{x^{2}}\partial_{x^{1}}p\left(x^{1}x^{2},x^{3}-x^{4}\right) & = & 1\\
\partial_{x^{3}}p\left(x^{1}x^{2},x^{3}-x^{4}\right) & = & 0
\end{eqnarray*}
and
\begin{eqnarray*}
\frac{1}{x^{2}}\partial_{x^{1}}p\left(x^{1}x^{2},x^{3}-x^{4}\right) & = & 0\\
\partial_{x^{3}}p\left(x^{1}x^{2},x^{3}-x^{4}\right) & = & 1
\end{eqnarray*}
we derive the flat outputs
\begin{eqnarray*}
y^{1} & = & x^{1}x^{2}\\
y^{2} & = & x^{3}-x^{4}\;.
\end{eqnarray*}
Finally it is worth mentioning that all distributions of both examples
are adjusted, to simplify the calculations. To determine the symmetries
in general, we need the functions, which generate $R^{\perp}$.

\section{Summary}

In this contribution we have presented two typical representations
of nonlinear sampled data systems, which turn out to be crucial for
their flatness. The main observation is, that flat systems must admit
a series or partial series connection of a Brunovsky form and a complement.
This fact follows from an integrability condition in the space of
the flat outputs and their shifts. Since the series connections admit
a simplification by reduction of the original problem to a simple
one, we get an algorithm to derive the flat outputs. A spin off of this
algorithm is the fact, that flat outputs are functions of the states.
The only exception is, if the original system has redundant input.
The disadvantage of the first version of the algorithm is, that one
has to solve linear PDEs or nonlinear ODEs, but one gets the flat
outputs directly. Therefore, a second version has been presented,
where at least a test for flatness is possible without the disadvantage
of the first one.


\bibliographystyle{plain}        
\bibliography{autosam}           



\end{document}